\documentclass [a4paper,twoside,11pt]{article}

\usepackage[latin1]{inputenc}
\usepackage{amsmath,amsfonts,amssymb}
\usepackage{vmargin,graphicx,theorem}
\usepackage[english]{babel}
\usepackage{enumerate}

\setpapersize[portrait]{A4}
\setmarginsrb{2.2cm}{2cm}{2.2cm}{1cm}{0cm}{0.1cm}{0.5cm}{2cm}

\newcommand{\disp}{\displaystyle}

\newcommand{\dN}{\ensuremath{\mathbb{N}}}

\newcommand{\dR}{\ensuremath{\mathbb{R}}}

\newtheorem{ethm}{Theorem}[section]

\newtheorem{ecor}[ethm]{Corollary}

\newtheorem{eprop}[ethm]{Proposition}

\newtheorem{elem}[ethm]{Lemma}

\newtheorem{erem}[ethm]{Remark}

\newcommand{\proofend}{~$\rhd$}
\newcommand{\proofbegin}{~$\lhd$}

\newenvironment{eproof}
               {\noindent {\emph{\textbf{Proof}}}\\\proofbegin~}
               {\proofend\\}

\newcommand{\p}[4]{{#3}\!\left#1{#4}\right#2}

\newcommand{\ABS}[1]{\ensuremath{{\left| #1 \right|}}} 
\newcommand{\PAR}[1]{\ensuremath{{\left(#1\right)}}} 
\newcommand{\BRA}[1]{\ensuremath{{\left\{#1\right\}}}} 
\newcommand{\NRM}[1]{\ensuremath{{\left\Vert #1\right\Vert}}} 

\renewcommand{\phi}{\varphi}

\renewcommand{\geq}{\geqslant}

\newcommand{\varf}[1]{\mathbf{Var}_{#1}}
\newcommand{\entf}[1]{\mathbf{Ent}_{#1}}

\newcommand{\ent}[2]{\p(){\entf{#1}}{#2}}

\newcommand{\var}[2]{\p(){\varf{#1}}{#2}}

\def\disp{\displaystyle}

\newcommand{\Ga}{\boldsymbol{\Gamma}}
\newcommand{\gu}{\boldsymbol{\Ga}}
\newcommand{\gd}{{\bf \gu {\!\!_2}}}

\newcommand{\GU}{\mathbf{\Gamma}}
\newcommand{\GD}{\GU \!\!_{\mathbf{2}}}

\newcommand{\al}{\alpha}

\newcommand{\ep}{\epsilon}

\newcommand{\la}{\lambda}

\newcommand{\e}{\varepsilon}

\newcommand{\affil}[1]{{\small\sl #1}}
\newcommand{\email}[1]{{\small E-mail: {\textsf {#1}}}}
\newcommand{\http}[1]{{\small Internet: {\textsf {#1}}}}

\begin{document}

\title{\sl From the Pr\'ekopa-Leindler inequality to modified  logarithmic Sobolev inequality}
\author{
  Ivan Gentil\\
\affil{Ceremade (UMR CNRS no. 7534), Universit\'e Paris IX-Dauphine,}\\
\affil{Place de Lattre de Tassigny, 75775 Paris C\'edex~16, France}\\
\email{gentil@ceremade.dauphine.fr}\\
\http{http://www.ceremade.dauphine.fr/
\raisebox{-4pt}{$\!\!\widetilde{\phantom{x}}$}gentil/}\\
}
\date{\today}\maketitle\thispagestyle{empty}

\begin{abstract}
We develop in this paper an improvement of the method given by S.
Bobkov and M. Ledoux in~\cite{bobkov-ledoux1}. Using the
Pr\'ekopa-Leindler inequality, we prove a modified logarithmic Sobolev
inequality adapted for all  measures  on $\dR^n$, with a strictly convex and super-linear potential. This
inequality implies modified logarithmic Sobolev inequality, developed in~\cite{ge-gu-mi,ge-gu-mi2},
for all uniformly strictly convex
potential as well as the  Euclidean logarithmic Sobolev inequality. 

\begin{center}
{\bf Résumé}
\end{center}

Dans cet article nous  am\'elirons la m\'ethode expos\'ee par S. Bobkov et M. Ledoux dans~\cite{bobkov-ledoux1}. En utilisant l'in\'egalit\'e de Pr\'ekopa-Leindler, nous prouvons une in\'egalit\'e de  Sobolev logarithmique modifi\'ee, adapt\'ee \`a toutes les mesures sur $\dR^n$ poss\'edant un  potentiel strictement convexe et super-lin\'eaire. Cette in\'egalit\'e implique en particulier une in\'egalit\'e de Sobolev logarithmique modifi\'ee, d\'evelopp\'ee dans~\cite{ge-gu-mi,ge-gu-mi2},  pour les mesures ayant un potentiel uniform\'ement strictement convexe mais aussi une in\'egalit\'e de Sobolev logarithmique de type euclidien. 
\end{abstract}

\section{Introduction}
\label{sec-int}

The Pr\'ekopa-Leindler inequality is the functional form of Brunn-Minkowski
inequality. Let $a,b$ be some positive reals such that $a+b=1$, 
and $u$, $v$, $w$ be some  non-negative measurable functions on $\dR^n$. Assume that, for any
$x,y \in \dR ^n$, we have
\begin{equation*}
 u(x)^a v(y)^b\leq w(a x + b y),
\end{equation*}
then
\begin{equation}
\label{ul-6.5} \left(\int u(x) dx\right)^a\left(\int v(x)
dx\right)^b\leq \int w(x) dx,
\end{equation}
where $dx$ is the Lebesgue measure on $\dR^n$. If we apply  inequality~\eqref{ul-6.5} to characteristic
functions of bounded measurable sets $A$ and $B$ in $\dR^n$, we get 
the multiplicative form of the Brunn-Minkowski inequality
$$
vol(A)^a vol(B)^b\leq vol(aA+bB),
$$
where $aA+bB=\BRA{ax_A+bx_B,\,\, x_A\in A,x_B\in B}$ and $vol(A)$ is the Lebesgue measure of the set $A$.
 One can see
for example two interesting reviews on this topic \cite{gupta,maurey}.

Bobkov and Ledoux in \cite{bobkov-ledoux1} use the Pr\'ekopa-Leindler
inequality to prove some functional inequalities like Brascamp-Lieb,
Logarithmic Sobolev and Transportation inequalities.

\medskip
More precisely,  let $\phi$ be a $\mathcal C^2$ strictly convex
function on $\dR^n$ and let
\begin{equation}
\label{eq-defm}
d\mu_\phi(x)=e^{-\phi(x)}dx
\end{equation}
be a probability measure on $\dR^n$ ($\int
e^{-\phi(x)}dx=1$). The function $\phi$ is called {\it the potential} of the measure $\mu_\phi$.  Bobkov and 
Ledoux obtained in particular the following two results:
\begin{itemize}
\item (Proposition~2.1 of \cite{bobkov-ledoux1}) Brascamp-Lieb inequality: assume that $\phi$ is a $\mathcal C^2$ function
on $\dR^n$, then for all smooth enough functions $g$,
\begin{equation}
\label{eq-br} \var{\mu_\phi}{g}:=\int \PAR{g-\int
gd\mu_\phi}^2d\mu_\phi\leq\int \nabla
g\cdot{\text{Hess}}(\phi)^{-1}\nabla gd\mu_\phi,
\end{equation}
where $\text{Hess}(\phi)^{-1}$ is the inverse of the Hessian of
$\phi$.
\item (Proposition~3.2 of \cite{bobkov-ledoux1}) Assume that for some $c>0$ and $p\geq2$,
for all $t,s>0$ with $t+s=1$, and for all $x,y\in\dR^n$, $\phi$
satisfies, as $s$ goes to 0,
\begin{equation}
\label{e-in}
t\phi(x)+s\phi(y)-\phi(tx+sy)\geq\frac{c}{p}(s+o(s))\NRM{x-y}^p,
\end{equation}
where $\NRM{\cdot}$ is the Euclidean norm in $\dR^n$. Then for all
smooth enough functions~$g$,
\begin{equation}
\label{e-bl1} \ent{\mu_\phi}{e^g}:=\int e^g\log \frac{e^g}{\int
e^gd\mu_\phi}d\mu_\phi\leq c\int \NRM{\nabla g}^qe^gd\mu_\phi,
\end{equation}
where $1/p+1/q=1$. They also give an example:  the function $\phi(x)=\NRM{x}^p+Z_\phi$ (where $Z_\phi$ is 
a normalization constant) which
satisfies inequality~\eqref{e-in} for some constant $c>0$.
\end{itemize}

The main result of this paper is to prove an inequality satisfies for any measure $\mu_\phi$ with a potential strictly convex and super-linear (we also assume a technical hypothesis satisfied by the potential $\phi$)).  More precisely we obtain, for all smooth enough  functions $g$ on~$\dR^n$,
\begin{equation}
\label{e-theo0} \ent{\mu_\phi}{e^g}\leq \int \BRA{x\cdot\nabla
g(x)-\phi^*(\nabla\phi(x))+\phi^*\PAR{\nabla\phi(x)-\nabla
g(x)}}e^{g(x)}d\mu_\phi(x),
\end{equation}
where $\phi^*$ is the Fenchel-Legendre transform of $\phi$,
$\phi^*(x):=\sup_{z\in\dR^n}\BRA{x\cdot z-\phi(z)}$.

The main application of this result is to extend the modified logarithmic Sobolev inequalities presented in 
\cite{ge-gu-mi,ge-gu-mi2} for probability measures on $\dR$ satisfying a
 uniform strictly convexity condition. It is well
 known that if the potential $\phi$ is $\mathcal C^2$ on $\dR$ such that 
for all $x\in\dR$, $\phi''(x)\geq\la>0$, then   the measure $\mu_\phi$ defined 
on~\eqref{eq-defm} verifies the logarithmic Sobolev inequality introduced by Gross in~\cite{gross}, for all 
smooth enough functions $g$, namely 
\begin{equation*}
\ent{\mu_\phi}{e^g}\leq \frac{1}{2\la}\int g'^2e^gd\mu_\phi.
\end{equation*}
This result comes from the $\GD$-criterion of D. Bakry and M. 
\'Emery, see \cite{b-e} or \cite{logsob} for a review.
We then improve the classical  logarithmic Sobolev inequality of Gross, in the situation where 
if the potential is even with $\phi(0)=0$ and satisfies 
\begin{equation*}
\forall x\in\dR, \quad \phi''(x)\geq \lambda >0\text{ and } \lim_{\ABS{x}\rightarrow\infty} \phi''(x)=\infty.
\end{equation*}
Adding a technical hypothesis (see Section~\ref{sec-mod}), we show that for all smooth functions $g$,
\begin{equation*}
\ent{\mu_\phi}{e^g}\leq \int H_\phi(g')e^gd\mu_\phi,
\end{equation*}
where 
\begin{equation*}
H_{\phi}(x)=
\left\{
\begin{array}{r}
\displaystyle C'\phi^*\PAR{\frac{x}{2}},\quad \text{if}\quad\ABS{x}> C\\
\displaystyle\frac{1}{2\la}x^2,\quad \text{if}\quad\ABS{x}\leq C,\\
\end{array}
\right.
\end{equation*}
for some constants $C,C',\la>0$ depending on  $\phi$. Remark that we always have 
$$
\forall x\in\dR,\quad
H_{\phi}(x)\leq C''x^2,
$$
for some other constant $C''$.
This inequality implies concentration inequalities which are more adapted 
to the measure studied, as we will see in Section~\ref{sec-mod}.

\bigskip 

The next section is divided into two subsections. In the first one we state the main theorem of this article, inequality~\eqref{e-theo0}. In the second subsection, we explain how this result improves results of \cite{bobkov-ledoux1}. In particular, inequality~\eqref{e-bl1} or Brascamp-Lieb inequality~\eqref{eq-br}. Section~\ref{sec-app} deals with some applications. The first one is an improvement of the a classical consequence of the $\gd$-criterion of Bakry-\'Emery for measures on $\dR$. We obtain then a global view of modified logarithmic Sobolev inequality for log-concave measures as introduced in joint work with A. Guillin and L. Miclo in \cite{ge-gu-mi,ge-gu-mi2}. Finally, we explain how the main theorem is equivalent to the Euclidean logarithmic Sobolev inequality. As a consequence, a short proof of the generalization given in \cite{del-dol,gentil03,ghoussoub} is obtained.

\section{Inequality for log-concave measures}
\label{sec-theo}

\subsection{The main theorem}
\label{s-re}

\begin{ethm}
\label{the-theo}
Let $\phi$ be a $\mathcal C^2$ strictly convex function on~$\dR^n$, such that
\begin{equation}
\label{eq-propphi}
\lim_{\NRM{x}\rightarrow\infty}\frac{\phi(x)}{\NRM{x}}=\infty.
\end{equation}
Denotes by $\mu_\phi(dx)=e^{-\phi(x)}dx$ a probability measure on~$\dR^n$, 
where $dx$ is the Lebesgue measure on~$\dR^n$, ($\int e^{-\phi(x)}dx=1$). Assume that $\mu_\phi$ satisfies for any $R>0$,
\begin{equation}
\label{hypo-++}
\int \PAR{\NRM{z}+\NRM{y_0}+R}^2\PAR{\NRM{\nabla\phi(z)}+\!\!\!\!\! \sup_{y;\,\NRM{y-z+y_0}\leq R} \!\!\!\!\!\NRM{\text{Hess}(\phi)(y)}}d\mu_\phi(z)<+\infty,
\end{equation}
where $y_0$ satisfies $\NRM{\nabla\phi(y_0)}\leq\NRM{\nabla\phi(z)}+R$.

If  
$\phi^*$ is the Fenchel-Legendre transform of $\phi$, 
$\phi^*(x):=\sup_{z\in\dR^n}\BRA{x\cdot z-\phi(z)},$ 
then for all smooth enough functions $g$ on $\dR^n$, one gets
\begin{equation}
\label{eq-theo} \ent{\mu_\phi}{e^g}\leq \int \BRA{x\cdot\nabla
g(x)-\phi^*(\nabla\phi(x))+\phi^*\PAR{\nabla\phi(x)-\nabla
g(x)}}e^{g(x)}d\mu_\phi(x).
\end{equation}
\end{ethm}

\begin{elem}
\label{lem-1}
Let $\phi$ satisfying conditions on Theorem~\ref{the-theo} then we have
\begin{itemize}
\item $\nabla \phi$ is a bijection on $\dR^n$ to $\dR^n$
\item $\lim_{\NRM{x}\rightarrow\infty}\frac{x\cdot\nabla\phi(x)}{\NRM{x}}=+\infty$.
\end{itemize}
\end{elem}

\begin{eproof}
Condition~\eqref{eq-propphi} implies that for all $x\in\dR^n$ the supremum of ${x\cdot z-\phi(z)}$ for $y\in\dR^n$ is reached for some $y\in\dR^n$.  Then $y$ satisfies $x=\nabla\phi(y)$ and  it proves that $\nabla\phi$ is a surjection. Then the strict convexity of $\phi$ implies that $\nabla\phi$ is a bijection. 

The function $\phi$ is convex then for all $x\in\dR^n$, $x\cdot\nabla\phi(x) \geq \phi(x)-\phi(0)$,~\eqref{eq-propphi} implies the second properties satisfied  by $\phi$.
\end{eproof}

The proof of the theorem is based on the following lemma:
\begin{elem}
\label{lem-technique}
Let $g$ be a $\mathcal C^\infty$ function with a compact support on $\dR^n$. Let $s,t\geq 0$ with $t+s=1$ and
denotes
$$
\forall z\in\dR^n,\quad g_s(z)=\sup_{z=tx+sy}\PAR{g(x)-\PAR{t\phi(x)+s\phi(y)-\phi(tx+sy)}}.
$$
Then there exists $R\geq 0$ such that, when $s$ goes to $0$, 
\begin{multline*}
 g_s(z)=g(z)+s\BRA{z\cdot\nabla
g(z)-\phi^*\PAR{\nabla\phi(z)}+\phi^*\PAR{\nabla\phi(x)-\nabla
g(x)}}\\
+ \PAR{(\NRM{z}+\NRM{y_0}+R)\NRM{\nabla \phi(z)}+(\NRM{z}+\NRM{y_0}+R)^2\!\!\!\!\!\sup_{y;\, \NRM{y-z+y_0}\leq R} \!\!\!\!\! \NRM{\text{Hess}(\phi)(y)}}{O(s^2)},
\end{multline*}
where $y_0$ satisfies $\NRM{\nabla\phi(y_0)}\leq\NRM{\nabla\phi(z)}+R$ and $O(s^2)$ is uniform on $z\in\dR^n$.
\end{elem}

\begin{eproof}
Let $s\in]0,1/2[$ and  $x=z/t-(s/t)y$, hence
$$
g_s(z)=\phi(z)+\sup_{y\in\dR^n}\PAR{g\PAR{\frac{z}{t}-\frac{s}{t}y}-t\phi\PAR{\frac{z}{t}-\frac{s}{t}y}-s\phi(y)}.
$$
Due to the fact that $g$ has a compact support and by  property~\eqref{eq-propphi} there exists $y_s\in\dR^n$ such that
$$
\sup_{y\in\dR^n}\PAR{g\PAR{\frac{z}{t}-\frac{s}{t}y}-t\phi\PAR{\frac{z}{t}-\frac{s}{t}y}-s\phi(y)}=
g\PAR{\frac{z}{t}-\frac{s}{t}y_s}-t\phi\PAR{\frac{z}{t}-\frac{s}{t}y_s}-s\phi(y_s).
$$
Moreover, $y_s$ satisfies
\begin{equation}
\label{eq-ys}
\nabla g\PAR{\frac{z}{t}-\frac{s}{t}y_s}-t\nabla\phi\PAR{\frac{z}{t}-\frac{s}{t}y_s}+t\nabla \phi(y_s)=0.
\end{equation}
Lemma~\ref{lem-1} implies that there exists a unique solution $y_0$ of the equation
\begin{equation}
\label{eq-der} \nabla\phi(y_0)={\nabla\phi(z)-\nabla g(z)},\quad
y_0=\PAR{\nabla\phi}^{-1}\PAR{\nabla\phi(z)-\nabla g(z)}.
\end{equation}

We prove now that $\lim_{s\rightarrow 0}y_s=y_0$.

\medskip

First we show  that there exists $A\geq0$ such that $\forall
s\in]0,1/2[$, $\NRM{y_s}\leq A$. Indeed, if the function $y_s$ is not bounded one can found $(s_k)_{k\in\dN}$ such that $s_k\rightarrow0$ and $\NRM{y_{s_k}}\rightarrow\infty$. 
Definition of $y_s$ implies that 
$$
g\PAR{\frac{z}{t}-\frac{s}{t}y_s}-t\phi\PAR{\frac{z}{t}-\frac{s}{t}y_s}-s\phi(y_s)\geq g\PAR{\frac{z}{t}}-t\phi\PAR{\frac{z}{t}}.
$$
Due to the fact that  $\lim_{\NRM{x}\rightarrow\infty}\phi(x)=\infty$ and since $g$ is bounded we obtain $s_ky_{s_k}=O(1)$. Next using~\eqref{eq-ys} one get
$$
\frac{y_s\cdot \nabla g\PAR{\frac{z}{t}-\frac{s}{t}y_s}}{\NRM{y_s}}-t\frac{y_s\cdot\nabla\phi\PAR{\frac{z}{t}-\frac{s}{t}y_s}}{\NRM{y_s}}+t\frac{y_s\cdot\nabla \phi(y_s)}{\NRM{y_s}}=0.
$$ 
The last equality is an contradiction with the second assertion of Lemma~\ref{lem-1} which prove that the $(y_s)_{s\in]0,1/2[}$ is bounded.

Let $\hat{y}$ an accumulation point of the function $y_s$, when  $s$ tends to 0.
Then $\hat{y}$ satisfies equation~\eqref{eq-der}. By unicity of
the solution of~\eqref{eq-der} we get  $\hat{y}=y_0$. Therefore we have
proved that $\lim_{s\rightarrow 0}y_s=y_0$.
\medskip

Taylor formula gives
$$
\phi\PAR{\frac{z}{t}-\frac{s}{t}y_s}=\phi(z)+s\PAR{\frac{z}{t}-\frac{y_s}{t}}\cdot\nabla\phi(z) + 
s^2\int_0^1(1-t)\PAR{\frac{z}{t}-\frac{y_s}{t}}\!\cdot\!{\rm{Hess}}(\phi)\PAR{\frac{z}{t}-s\frac{y_s}{t}}\PAR{\frac{z}{t}-\frac{y_s}{t}}dt,
$$
and the same for $g$. Using  the continuity of $y_s$ at $s=0$, ones gets
\begin{multline*}
\phi\PAR{\frac{z}{t}-\frac{s}{t}y_s}=\phi(z)+s(z-y_0)\cdot\nabla\phi(z)\\ 
+ \PAR{(z-y_0)\cdot\nabla\phi(z)+\sup_{t\in[0,1/2]}\NRM{\frac{z}{t}-\frac{y_t}{t}}^2\sup_{t\in[0,1/2]}\NRM{{\rm{Hess}}(\phi)\PAR{\frac{z}{t}-\frac{y_s}{t}}}}{O(s^2)}.
\end{multline*}
and the same for $g$
\begin{multline*}
g\PAR{\frac{z}{t}-\frac{s}{t}y_s}=g(z)+s(z-y_0)\cdot\nabla g(z)
+\\
\PAR{(z-y_0)\cdot\nabla g(z)+\sup_{t\in[0,1/2]}\NRM{\frac{z}{t}-\frac{y_t}{t}}^2\sup_{t\in[0,1/2]}\NRM{{\rm{Hess}}(g)\PAR{\frac{z}{t}-\frac{y_s}{t}}}}{O(s^2)}.
\end{multline*}
As a consequence, 
\begin{multline*}
g_s(z)=g(z)+s\BRA{\phi(z)-\phi(y_0)+(z-y_0)\cdot(\nabla
g(z)-\nabla\phi(z))}\\
+ \PAR{(z-y_0)\cdot(\nabla g(z)-\nabla\phi(z))+\sup_{t\in[0,1/2]}\NRM{\frac{z}{t}-\frac{y_t}{t}}^2\sup_{t\in[0,1/2]}\NRM{{\rm{Hess}}(\phi+g)\PAR{\frac{z}{t}-\frac{y_s}{t}}}}{O(s^2)}.
\end{multline*}
The function $g$ is $\mathcal C^\infty$  with a compact support then one obtains using~\eqref{eq-der} and the expression of the
Fenchel-Legendre transformation for a strictly convex function
$$
\forall x\in\dR^n,\quad\phi^*(\nabla \phi (z))=\nabla\phi(z)\cdot z-\phi(z),
$$
we get the result.
\end{eproof}

We are now ready to deduce our main result:

{\noindent {\emph{\textbf{Proof of
Theorem~\ref{the-theo}}}\\\proofbegin~} The proof is based on the
proof of Theorem~3.2 of \cite{bobkov-ledoux1}. First we prove
inequality~\eqref{eq-theo} for all functions $g$, $\mathcal
C^\infty$ with  compact support on $\dR^n$.

 Let $t,s\geq 0$ with $t+s=1$ and denote for $z\in\dR^n$,
$$
g_t(z)=\sup_{z=tx+sy}\PAR{g(x)-\PAR{t\phi(x)+s\phi(y)-\phi(tx+sy)}}.
$$
We apply Pr\'ekopa-Leindler inequality to the functions
$$
u(x)=\exp\PAR{\frac{g(x)}{t}-\phi(x)},\quad v(y)=\exp\PAR{-\phi(y)},\quad w(z)=\exp\PAR{g_s(z)-\phi(z)},
$$
to get
$$
\PAR{\int \exp(g/t)d\mu_\phi}^t\leq\int \exp(g_s)d\mu_\phi.
$$
The differentiation of the $L^p$ norm gives the entropy, and thanks to a 
Taylor's formula we get
$$
\PAR{\int \exp(g/t)d\mu_\phi}^t=\int e^gd\mu_\phi+s\ent{\mu_\phi}{e^g}+O(s^2).
$$
Then applying Lemma~\ref{lem-technique} and inequality~\eqref{hypo-++} yield
\begin{multline*}
\int \exp(g_s)d\mu_\phi=\\
\int e^g\mu_\phi+ s\int
\BRA{z\cdot\nabla g(z)-\phi^*\PAR{\nabla\phi(z)}+
\phi^*\PAR{\nabla\phi(z)-\nabla g(z)}}e^{g(z)} d\mu_\phi(z)+O(s^2).
\end{multline*}
When $s$ goes to 0, inequality~\eqref{eq-theo} arises and can be extended for all  smooth enough functions $g$. {\proofend\\}

Note that hypothesis~\eqref{hypo-++} is satisfied by a large class of convex functions. For example if $\phi(x)=\NRM{x}^2/2+(n/2)\log(2\pi)$ we obtain the
classical logarithmic Sobolev of Gross for the canonical Gaussian
measure on $\dR^n$, with the optimal constant.

\subsection{Remarks and examples} \label{sec-ex}

In the  next corollary we recall a classical result of
perturbation.
If $\Phi$ is a function on $\dR^n$ such that $\int e^{-\Phi} dx<\infty$ we note the probability
 measure $\mu_\Phi$  by 
\begin{equation}
\label{e-phi}
d\mu_\Phi(x)=\frac{e^{-\Phi(x)}}{Z_\Phi}dx,
\end{equation}
where $Z_\Phi=\int {e^{-\Phi(x)}}dx$·

\begin{ecor}
\label{co-pe} Assume that $\phi$ satisfies conditions of Theorem~\ref{the-theo}. 
Let $\Phi=\phi+U$, where $U$ is a bounded function on $\dR^n$ and
denote by $\mu_\Phi$ the measure defined by~\eqref{e-phi}.

Then for all smooth enough functions $g$ on $\dR^n$, one has 
\begin{equation}
\label{eq-core} \ent{\mu_\Phi}{e^g}\leq e^{2\text{osc}(U)}\int
\BRA{x\cdot\nabla
g(x)-\phi^*(\nabla\phi(x))+\phi^*\PAR{\nabla\phi(x)-\nabla
g(x)}}e^{g(x)}d\mu_\Phi(x),
\end{equation}
where $\text{osc}(U)=\sup(U)-\inf(U)$.
\end{ecor}

\begin{eproof}
First we observe that
\begin{equation}
\label{eq-ra}
e^{-\text{osc}(U)}\leq\frac{d{\mu}_\Phi}{d\mu_\phi}\leq
e^{\text{osc}(U)}.
\end{equation}
Moreover we have  for all probability measures $\nu$ on $\dR^n$,
$$
\ent{\nu}{e^g}=\inf_{a\geq0}\BRA{\int\PAR{e^g\log\frac{e^g}{a}-e^g+a}d\nu}.
$$
Using the fact that for all $x,a>0$, $x \log\frac{x}{a}-x+a\geq0,$
we get
$$
e^{-\text{osc}(U)}\ent{{\mu}_\Phi}{e^g}\leq\ent{\mu_\phi}{e^g}\leq
e^{\text{osc}(U)}\ent{{\mu}_\Phi}{e^g}.
$$
Then if $g$ a smooth enough function  on $\dR^n$ we have
\begin{equation*}
\begin{array}{rl}
\disp\ent{{\mu}_\Phi}{e^g}& \disp\leq e^{\text{osc}(U)}\ent{{\mu}_\phi}{e^g}\\
&\disp \leq  e^{\text{osc}(U)} \int \BRA{x\cdot\nabla
g(x)-\phi^*(\nabla\phi(x))+\phi^*\PAR{\nabla\phi(x)-\nabla
g(x)}}e^{g(x)}d\mu_\phi(x).
\end{array}
\end{equation*}
The convexity of $\phi^*$  $\dR^n$ and th  relation
$\nabla\phi^*\PAR{\nabla\phi(x)}=x$ lead to
$$
\forall x\in\dR^n,\quad {x\cdot\nabla
g(x)-\phi^*(\nabla\phi(x))+\phi^*\PAR{\nabla\phi(x)-\nabla
g(x)}}\geq0.
$$
Finally by~\eqref{eq-ra} we get
\begin{equation*}
\ent{{\mu}_\Phi}{e^g} \leq  e^{2\text{osc}(U)} \int \BRA{x\cdot\nabla
g(x)-\phi^*(\nabla\phi(x))+\phi^*\PAR{\nabla\phi(x)-\nabla
g(x)}}e^gd\mu_\Phi.
\end{equation*}
\end{eproof}

\begin{erem}
It is not necessary to state a tensorization result, as we may
obtain exactly the same expression when  computing directly with a
product measure.
\end{erem}

Theorem~\ref{the-theo} implies also examples given
in \cite{bobkov-ledoux1} and \cite{bobkov-zeg}.
\begin{ecor}[\cite{bobkov-ledoux1}]
\label{bl-prop} Let $p\geq 2$ and let $\Phi(x)=\NRM{x}^p/p$ where
$\NRM{\cdot}$ is Euclidean norm in $\dR^n$. Then there exists $c>0$, such that  for all
smooth enough functions $g$,
\begin{equation}
\label{e-prop} \ent{\mu_\Phi}{e^g}\leq c\int \NRM{\nabla
g}^qe^gd\mu_\Phi,
\end{equation}
where $1/p+1/q=1$  and $\mu_\Phi$ is defined on~\eqref{e-phi}. 
\end{ecor}

\begin{eproof}
Using Theorem~\ref{the-theo}, we just have  to prove that there exists $c>0$ such that, 
\begin{equation*}
\forall x,y\in\dR^n,\quad {x\cdot y-\Phi^*(\nabla\Phi(x))+\Phi^*\PAR{\nabla\Phi(x)-
y}}\leq c\NRM{y}^q.
\end{equation*}
Assume that $y\neq 0$ and define the function $\psi$ by,
$$
\psi(x,y)=\frac{{x\cdot y-\Phi^*(\nabla\Phi(x))+\Phi^*\PAR{\nabla\Phi(x)-y}}}{\NRM{y}^q}.
$$
Then $\psi $ is a bounded function. 
We know that $\Phi^*(x)=\NRM{x}^q/q$. Choosing
 $z=x\NRM{x}^{p-2}/\NRM{y}$ and denoting $e=y/\NRM{y}$, we obtain
$$
\psi(x,y)=\bar{\psi}(z,e)=z\cdot
e\NRM{z}^{q-2}-\frac{1}{q}\NRM{z}^q+\frac{1}{q}\NRM{{z}-{e}}^q.
$$
Taylor's formula then yields 
$\bar{\psi}(z,e)=O(\NRM{z}^{q-2})$. But $p\geq 2$ implies that
$q\leq2$, so that $\bar{\psi}$ is a bounded function. We then  get 
the result with $c=\sup\bar{\psi}=\sup{\psi}$.
\end{eproof}

Optimal transportation  is also used by Cordero-Erausquin, Gangbo and Houdr\'e in~\cite{co-ga-ho} to  prove the particular case of
the inequality~\eqref{e-prop}.
\medskip

In  Proposition~2.1 of \cite{bobkov-ledoux1}, Bobkov and Ledoux
prove that  the Pr\'ekopa-Leindler inequality implies Brascamp-Lieb
inequality. In our case, Theorem~\ref{the-theo} also
implies   Brascamp-Lieb inequality, as we can see in the
next corollary.

\begin{ecor}
Let $\phi$ satisfying  conditions of Theorem~\ref{the-theo}.  Then
for all smooth enough functions $g$ we get,
\begin{equation*}
\var{\mu_\phi}{g}\leq \int \nabla
g\cdot{\text{Hess}}(\phi)^{-1}\nabla gd\mu_\phi,
\end{equation*}
where ${\text{Hess}}(\phi)^{-1}$ denote the inverse of the Hessian
of $\phi$.
\end{ecor}

\begin{eproof}
Assume that $g$ is a ${\mathcal C}^\infty$  function with a
compact support  and apply inequality~\eqref{eq-theo} with the
function $\ep g$ where $\ep>0$. Taylor's formula gives
$$
\ent{\mu_\phi}{\exp{\ep g}}=\frac{\ep^2}{2}\var{\mu_\phi}{g}+o(\ep^2),
$$
and
\begin{multline*}
\int \BRA{x\cdot\nabla
g(x)-\phi^*(\nabla\phi(x))+\phi^*\PAR{\nabla\phi(x)-\nabla
g(x)}}e^g{(x)}d\mu_\phi(x) =
\\\int\frac{\ep^2}{2}\nabla
g\cdot{\text{Hess}}(\phi^*)\PAR{\nabla \phi}\nabla
gd\mu_\phi+o(\ep^2).
\end{multline*}
Because of $\nabla\phi^*(\nabla\phi(x))=x$, one has 
${\text{Hess}}(\phi^*)\PAR{\nabla\phi}={\text{Hess}}(\phi)^{-1}$ which finished the proof. 
\end{eproof}

\begin{erem}
Let $\phi$ satisfying the conditions  of Theorem~\ref{the-theo}, and  $L$ be defined by
$$
\forall x,y\in\dR^n,\quad L(x,y)=\phi(y)-\phi(x)+(y-x)\cdot\nabla\phi(x).
$$
The convexity of $\phi$ implies  that $L(x,y)\geq0$ for all
$x,y\in\dR^n$. Let $F$ be a density of probability with respect to the measure
$\mu_\phi$, we defined the following Wasserstein distance with the
cost function $L$ by
$$
 W_L(Fd\mu_\phi,d\mu_\phi)=\inf\BRA{ \int L(x,y)d\pi(x,y)},
$$
where the infimum is taken over all probability measures $\pi$ on
$\dR^n\times\dR^n$ with marginal distributions $Fd\mu_\phi$ and
$d\mu_\phi$. Bobkov and Ledoux proved in \cite{bobkov-ledoux1}
the following transportation inequality
\begin{equation}
\label{eq-t} W_L(Fd\mu_\phi,d\mu_\phi)\leq \ent{\mu_\phi}{F}.
\end{equation}

The main result of Otto and Villani in \cite{villani} is the
following: Classical logarithmic Sobolev inequality of Gross  (when
$\phi(x)=\NRM{x}^2/2+(n/2)\log (2\pi)$) implies the transportation
inequality~\eqref{eq-t} for all functions $F$, density of
probability with respect to $\mu_\phi$ (see also \cite{bgl} for an another
 proof). The method developed in~\cite{bgl}, enables to
extend the property  for $\phi(x)=\NRM{x}^p+Z_\phi$ ($p\geq2$).
In the general case,   inequality proved in this article, we do not know if the modified logarithmic Sobolev 
inequality~\eqref{eq-theo} implies transportation inequality~\eqref{eq-t}.
\end{erem}

\section{Applications}
\label{sec-app}

\subsection{Application to modified logarithmic Sobolev inequalities}
\label{sec-mod}

In \cite{ge-gu-mi,ge-gu-mi2},  a modified logarithmic Sobolev inequality for measure $\mu_\phi$ on $\dR$  is given with a potential between $\ABS{x}$ and $x^2$. More precisely let $\Phi$ be a function on the real line and assume that $\Phi$ is even and satisfies the following property: there exist $M\geq 0$ and $0<\e\leq1/2$ such that,
\begin{equation*}
\tag{\bf H} \forall{x\geq M},\,\,\,\,(1+\e)\Phi(x)\leq
x\Phi'(x)\leq(2-\e)\Phi(x).
\end{equation*}

Then there exist $A,B,D>0$ such that for all smooth functions $g$
 we have
\begin{equation}
\label{ggl}
 \ent{\mu_\Phi}{e^g}\leq
A\int H_{\Phi}\PAR{g'}e^gd\mu_\Phi,
\end{equation}
where
\begin{equation*}
\label{defh}
H_{\Phi}(x)=\left\{
\begin{array}{rl}
\Phi^*\PAR{Bx} &\text{ if }\ABS{x}\geq D,\\
x^2 &\text{ if }\ABS{x}\leq D,
\end{array}
\right.
\end{equation*}
and $\mu_\Phi$ is defined on~\eqref{e-phi}.

The proof of inequality~\eqref{ggl} is rather technical and is
divided in two parts: the large and the small entropy. 
Using Theorem~\ref{the-theo} one obtains two results in this direction. In the next theorem,
 we extend~\eqref{ggl} in the  case where the potential 
is ``bigger'' than $x^2$.

\begin{ethm}
\label{theo-logg}
Let $\phi$ be a real function  satisfying conditions of Theorem~\ref{the-theo}. Assume that $\phi$ is even, $\phi(0)=0$,  $\phi''$ in decreasing on $]-\infty,0]$    and increasing on $[0,+\infty[$  and satisfies, 
\begin{equation}
\label{eq-hypo}
\forall x\in\dR, \quad \phi''(x)\geq \phi''(0)=\lambda>0 \text{ and }  \lim_{\ABS{x}\rightarrow\infty} \phi''(x)=\infty.
\end{equation}
Assume also   that there exists $A>1$ such that for $\ABS{x}\geq C$ for some $C>0$,
\begin{equation}
\label{eq-fa}
A\phi(x)\leq  x\phi'(x).
\end{equation}

Then there exists $C>0$ such that for all smooth enough functions $g$, 
\begin{equation}
\label{eq-logg}
\ent{\mu_\phi}{e^g}\leq \int H_\phi(g')e^gd\mu_\phi,
\end{equation}
where 
\begin{equation}
\label{eq-defhh}
H_{\phi}(x)=
\left\{
\begin{array}{r}
\displaystyle\frac{2A}{A-1}\phi^*\PAR{\frac{x}{2}},\quad \text{if}\quad\ABS{x}> C\\
\displaystyle\frac{1}{2\la}x^2,\quad \text{if}\quad\ABS{x}\leq C.\\
\end{array}
\right.
%
\end{equation}
\end{ethm}

The proof of this theorem is a straightforward   application of the following  lemma:
\begin{elem}
\label{lem-gg}
Assume that $\phi$ satisfies conditions of Theorem~\ref{theo-logg}, then we get
\begin{equation}
\label{eq-lemm}
\forall x,y\in\dR,\quad xy-\phi^*(\phi'(x))+\phi^*(\phi'(x)-y)\leq H_\phi(y).
\end{equation}
\end{elem}

\begin{eproof}
We know that for all $x\in\dR^n$, $x=\phi^{*'}(\phi'(x))$, and  the convexity of $\phi^*$ yields
\begin{equation}
\label{eq-bas}
xy-\phi^*(\phi'(x))+\phi^*(\phi'(x)-y)\leq y\PAR{\phi^{*'}(\phi'(x))-\phi^{*'}(\phi'(x)-y)}. 
\end{equation}
Let $y\in\dR$ be fixed and notes $\psi_y(x)=\phi^{*'}(x+y)-\phi^{*'}(x)$. 
The function $\phi$ is convex, so  one gets
for all $x\in\dR$, $\phi^{*''}(\phi'(x))\phi''(x)=1$, and
 the maximum of $\psi_y(x)$ is reached on $x_0\in\dR$ 
which satisfies the condition ${\phi^{*''}(x_0)=\phi^{*''}(x_0+y)}$. Since $\phi$ is even, $\phi''$ is decreasing on $]-\infty,0]$ and increasing on $]-\infty,0]$ one get that $x_0+y=-x_0$. Then one obtains 
$$
\forall y\in\dR,\quad\forall x\in\dR,\quad\phi^{*'}(x+y)-\phi^{*'}(x)\leq2\phi^{*'}\PAR{\frac{y}{2}},
$$


$$
\forall y\in\dR,\quad\forall x\in\dR,\quad xy-\phi^*(\phi'(x))+\phi^*(\phi'(x)-y)\leq y2\phi^{*'}\PAR{\frac{y}{2}}. 
$$ 
By \eqref{eq-fa} one gets 
\begin{equation}
\label{ea-dr}
\forall \ABS{y}\geq C,\quad\forall x\in\dR,\quad xy-\phi^*(\phi'(x))+\phi^*(\phi'(x)-y)\leq \frac{2A}{A-1}\phi^{*}\PAR{\frac{y}{2}}. 
\end{equation}  
A Taylor's formula then leads to 
$$
\forall x,y\in\dR,\quad xy-\phi^*(\phi'(x))+\phi^*(\phi'(x)-y)\leq \frac{y^2}{2}\phi^{*''}(\phi'(x)-\theta y)
\leq \frac{y^2}{2\la},
$$ 
for some $\theta\in(0,1)$, and one gets~\eqref{ea-dr}. 
\end{eproof}

\begin{erem}
\begin{itemize}
\item The last theorem improved the classical consequence of Bakry-\'Emery criterion for the logarithmic Sobolev inequality. In fact when a probability measure is more log-concave than the Gaussian measure, we obtain a modified logarithmic Sobolev inequality sharper than the classical inequality of Gross. Using a such  inequality then one obtains concentration inequality which is  more adapted to the probability measure studied. 
\item Theorem~\ref{theo-logg} is more precise than  Corollary~\ref{bl-prop}  proved by Bobkov, Ledoux and Zegarlinski in \cite{bobkov-ledoux1, bobkov-zeg}. The particularity of the function $H_\phi$ defined on~\eqref{eq-defhh} is its behaviour  around the origin. One can obtain easily that if a probability measure satisfies inequality~\eqref{eq-logg} then it satisfies  a Poincar\'e inequality with constant $1/\la$. 
\item Note also that this method can not be applied for measures with a concentration between $e^{-\ABS{x}}$ and $e^{-x^2}$ described in \cite{ge-gu-mi,ge-gu-mi2}.  In particular  Lemma~\ref{lem-gg} is false in this case. 
\item Note  finally that the condition~\eqref{eq-fa} is a technical condition, satisfied for a large class of functions.
\end{itemize}
\end{erem}

A natural application of Theorem~\ref{theo-logg} is a concentration inequality in the spirit of Talagrand, see~\cite{talagrand}.
\begin{ecor}
\label{cor-cons}
  Assume that $\phi$ satisfies  conditions of Theorem~\ref{theo-logg} and there exists $B>1$ such that for $\ABS{x}$ large enough,
\begin{equation}
\label{eq-fa2}
 x\phi'(x)\leq B\phi(x).
\end{equation}
Then there exists
constants $C_1,C_2,C_3\geq0$, independent of $n$ such that:
if $F$ is a function on $\dR^n$ such that $\forall i$, 
$\NRM{\partial_iF}_\infty\leq 1$, then we get for $\la\ge0$,
\begin{equation}
\label{eq-cons}
\mu^{\otimes n}(\ABS{F -\mu^{\otimes n}(F)}\geq \la)\leq
\left\{
\begin{array}{ll}
\disp2\exp\PAR{-nC_1\Phi\PAR{C_2\frac{\la}{n}}}&\text{if }
\la > {nC_3 },\\
\disp2\exp\PAR{-C_1\frac{\la^2}{n}}& \text{if }0
\leq\la\leq {nC_3 }.
\end{array}
\right.
\end{equation}
\end{ecor}

\begin{eproof}
Using the additional hypothesis~\eqref{eq-fa2}, the proof of~\eqref{eq-cons} is the same as for Proposition~3.2 of~\cite{ge-gu-mi2}.
\end{eproof}

A $n$-dimensional version of~\eqref{eq-logg} is also available. 
\begin{eprop}
\label{thmn} Let $\Phi$ be a $\mathcal C^2$,  strictly convex and
even function on $\dR^n$ and satisfying~\eqref{eq-propphi} and~\eqref{hypo-++}.
Assume also that $\Phi\geq0$ and $\Phi(0)=0$ {$($}it implies that $0$ is the unique
minimum of $\Phi${$)$},
\begin{equation}
\label{eqm}
\lim_{\al\rightarrow0,\,\al\in[0,1]}
\sup_{x\in\dR^n}\BRA{(1-\al)\frac{\Phi^*\PAR{\frac{x}{1-\al}}}{\Phi^*(x)}}=1,
\end{equation}
and also that there exists $A>0$ such that
\begin{equation}
\label{eqh}
\forall x\in\dR^n,\quad x\cdot \nabla \Phi(x)\leq (A+1)\Phi(x).
\end{equation}

Then there exist $C_1,C_2,C_3\geq0$ such that for all smooth enough functions $g$ such that
$\int e^gd\mu_\Phi=1$we get
\begin{equation}
\label{eqd}
\ent{\mu_\Phi}{e^g}\leq C_1\int\Phi^*\PAR{C_2\nabla g}e^gd\mu_\Phi+C_3.
\end{equation}
\end{eprop}

\begin{eproof}
Let apply Theorem~\ref{the-theo} with $\phi=\Phi+\log Z_\Phi$, one has
\begin{equation*}
\ent{\mu_\Phi}{e^g}\leq \int \BRA{x\cdot\nabla
g(x)-\Phi^*(\nabla\Phi(x))+\Phi^*\PAR{\nabla\Phi(x)-\nabla
g(x)}}e^gd\mu_\Phi.
\end{equation*}

The convexity of $\Phi^*$ implies, for all $\al\in[0,1[$,
\begin{equation}
\label{eqc}
\forall x\in\dR^n,\quad\Phi^*\PAR{\nabla\Phi(x)-\nabla g(x)}\leq(1-\al)\Phi^*\PAR{\frac{\nabla\Phi(x)}{1-\al}}+
\al\Phi^*\PAR{\frac{-\nabla g(x)}{\al}}.
\end{equation}
Recall that $\Phi^*$ is also an even function.
Young's inequality implies that
\begin{equation}
\label{eqy}
\forall x\in\dR^n,\quad x\cdot\frac{\nabla g(x)}{\al}\leq\Phi(x)+\Phi^*\PAR{\frac{\nabla g(x)}{\al}}.
\end{equation}
Using~\eqref{eqc} and ~\eqref{eqy} we get
\begin{multline*}
\ent{\mu_\Phi}{e^g}\leq 2\al\int \Phi^*\PAR{\frac{\nabla g}{\al}}e^gd\mu_\Phi+
\al \int \Phi\PAR{x}e^gd\mu_\Phi+\\
\int \PAR{(1-\al)\Phi^*\PAR{\frac{\nabla\Phi(x)}{1-\al}}-\Phi^*\PAR{{\nabla\Phi(x)}}}e^gd\mu_\Phi.
\end{multline*}
We have $\Phi^*(\nabla \Phi(x))=x\cdot\nabla \Phi(x)-\Phi(x)$, then
inequality~\eqref{eqh} implies that $\Phi^*\PAR{\nabla \Phi(x)}\leq A\Phi(x)$. Because of 
$\Phi(0)=0$ one has  $\Phi^*\geq 0$, so that 
$$
\ent{\mu_\Phi}{e^g}\leq \al\int \Phi^*\PAR{\frac{\nabla g}{\al}}e^gd\mu_\Phi+
 \al\int \Phi^*\PAR{\frac{\nabla g}{\al}}e^gd\mu_\Phi+
(\al+A\ABS{\psi(\al)-1})\int \Phi e^gd\mu_\Phi,
$$
where
\begin{equation}
\label{eqp}
\psi(\al)=\sup_{x\in\dR^n}\BRA{(1-\al)\frac{\Phi^*\PAR{\frac{x}{1-\al}}}{\Phi^*(x)}}.
\end{equation}
Let $\la>0$, recall that  $\int e^g d\mu_\Phi=1$ then yields 
$$
\int \Phi e^gd\mu_\Phi\leq \la\PAR{\ent{\mu_\Phi}{e^g}+\log\int e^{\Phi /\la}d\mu_\Phi}.
$$
One has $\disp\lim_{\la\rightarrow\infty}\log\int e^{\Phi /\la}d\mu_\Phi=0$. Let then 
let now choose $\la$ large enough so  that  $\log\int e^{\Phi /\la}d\mu_\Phi\leq 1$.
Using the property~\eqref{eqm}, taking $\al$ such that
$(\al+A\ABS{\psi(\al)-1})\la\leq1/2$ implies
$$
\ent{\mu_\Phi}{e^g}\leq 2\al\int \Phi^*\PAR{\frac{\nabla g}{\al}}e^gd\mu_\Phi+
\frac{1}{2}\PAR{\ent{\mu_\Phi}{e^g}+1},
$$
which gives
$$
\ent{\mu_\Phi}{e^g}\leq{4\al}\int \Phi^*\PAR{\frac{\nabla g}{\al}}e^gd\mu_\Phi+1/2.
$$
\end{eproof}

The main difference between the inequality obtained and the modified logarithmic inequality~\eqref{eq-logg} is that we do not have equality if $f=1$. Then~\eqref{eqd} is called a no tight inequality and it is more difficult to obtain. 

\subsection{Application to Euclidean logarithmic Sobolev inequality}
\label{sec-e}

\begin{ethm}
\label{theo2} Assume that the function $\phi$ satisfies conditions
of Theorem~\ref{the-theo}. Then for all $\la>0$ and for all smooth
enough functions $g$ on $\dR^n$,
\begin{equation}
\label{e-theo2}
\ent{dx}{e^g}\leq -n\log\PAR{\la e}\int e^gdx+\int\phi^*\PAR{-\la\nabla g}e^gdx.
\end{equation}

This inequality is optimal in the sense that if $g=-\phi(x-\bar{x})$
with $\bar{x}\in\dR^n$ and $\la=1$ we get an equality.
\end{ethm}

\begin{eproof}
Integrating by parts in the second term of~\eqref{eq-theo} yields for all $g$ smooth enough
$$
\int {x\cdot\nabla g(x)}e^{g(x)}d\mu_\phi(x)=
\int \PAR{-n+x\cdot\nabla\phi(x)}e^{g(x)}d\mu_\phi(x).
$$
Then using the equality $\phi^*\PAR{\nabla\phi}=x\cdot\nabla\phi(x)-\phi(x)$ we get
for all smooth enough functions $g$,
\begin{equation*}
\ent{\mu_\phi}{e^g}\leq
\int \PAR{-n+\phi+\phi^*\PAR{\nabla\phi-\nabla g}}e^{g}d\mu_\phi,
\end{equation*}
Let now take $g=f+\phi$ to obtain
\begin{equation*}
\ent{dx}{e^f}\leq
\int \PAR{-n+\phi^*\PAR{-\nabla g}}e^{g}dx.
\end{equation*}
Finally, let $\la>0$ and take $f(x)=g(\la x)$, we get then
\begin{equation*}
\ent{dx}{e^g}\leq
-n\log\PAR{\la e}\int e^gdx+\int\phi^*\PAR{-\la\nabla g}e^gdx,
\end{equation*}
which proves~\eqref{e-theo2}.

If now  $g=-\phi(x-\bar{x})$ with $\bar{x}\in\dR^n$ an easy computation proves that if $\la=1$ the 
equality holds.
 \end{eproof}

In the inequality~\eqref{e-theo2}, there exists an optimal
$\la_0>0$ and  when $C$ is homogeneous, we can improve the
last result. We find an inequality called {\it Euclidean logarithmic
Sobolev inequality} which is explained  on the next  corollary.

\begin{ecor}
\label{cor1}
Let $C$ be a strictly convex function on $\dR^n$ satisfying condition of Theorem~\ref{the-theo} and assume that $C$ is  $q$-homogeneous for some $q>1$,
$$
\forall \la\geq0\quad \text{and}\quad \forall x\in\dR^n,\quad C(\la x)={\la}^qC(x).
$$
Then for all smooth enough functions $g$ in $\dR^n$ we get
\begin{equation}
\label{e-ecl} \ent{dx}{e^g}\leq \frac{n}{p}\int
e^gdx\log\PAR{\frac{p}{n e^{p-1}{\mathcal L}^{p/n}}\frac{\int
C^*\PAR{-\nabla g}e^gdx}{\int e^gdx} },
\end{equation}
where ${\mathcal L}=\int e^{-C}dx$ and $1/p+1/q=1$.
\end{ecor}

\begin{eproof}
Let apply Theorem~\ref{theo2} with $ \phi=C+\log{\mathcal L}$. Then $\phi$ satisfies conditions
of Theorem~\ref{theo2} and we get then
\begin{equation*}
\ent{dx}{e^g}\leq -n\log\PAR{\la e {\mathcal L}^{1/n}}\int
e^gdx+\int C^*\PAR{-\la\nabla g}e^gdx.
\end{equation*}
Due to the fact that $C$ is $q$-homogeneous an easy computation proves
that  $C^*$ is  $p$-homogeneous where $1/p+1/q=1$. An optimization
over $\la>0$ gives inequality~\eqref{e-ecl}.
\end{eproof}

\begin{erem}
Inequality~\eqref{e-ecl} is useful to prove regularity properties as  hypercontractivity  for nonlinear diffusion as the $p$-Laplacian, see~\cite{ddg}.  The function $C$ is then adapted to the nonlinear diffusion studied.  
\end{erem}

Inequality~\eqref{e-ecl} is called Euclidean logarithmic Sobolev
inequality and computations of this section  is the generalization of the 
work of Carlen in~\cite{carlen}. This inequality with $p=2$, appears in the work of
Weissler in~\cite{weissler}.
It was discussed and extended to this
last version in many articles see
\cite{beckner,del-dol,gentil03,ghoussoub}.

\begin{erem}
As explained in the introduction, computation used in Corollary~\ref{cor1} clearly  proves that
inequality~\eqref{e-ecl} is equivalent to inequality~\eqref{e-theo2}.  Agueh, Ghoussoub and Kang, in~\cite{ghoussoub},
used Monge-Kantorovich theory for mass transport to prove inequalities~\eqref{e-theo2} and~\eqref{e-ecl}. Their approach
gives another way to establish Theorem~\ref{the-theo}.

Note finally  that inequality~\eqref{e-ecl} is optimal, extremal functions are given by
$g(x)=-bC(x-\bar{x})$, with $\bar{x}\in\dR^n$ and $b>0$. If they are only ones is still an open question.
\end{erem}

{\bf{Acknowledgments}}: 
I would like to warmly thank referee for pointed out errors in the first version. 

\newcommand{\etalchar}[1]{$^{#1}$}

\end{document}